# A Quick Response Algorithm for Dynamic Autonomous Mobile Robot Routing Problem with Time Windows

Lulu Cheng[a] Ning Zhao[a] Mengge Yuan[a] and Kan Wu[b]*

*[a]Faculty of Science, Kunming University of Science and Technology, Kunming, Yunnan;*
*[b]Business Analytics Research Center, Chang Gung University, Taoyuan City, Taiwan*
Corresponding author: Kan Wu; E-mail: kan626@gmail.com; Business Analytics Research Center, Chang Gung University, Taoyuan City, Taiwan

# A Quick Response Algorithm for Dynamic Autonomous Mobile Robot Routing Problem with Time Windows

This paper investigates the optimization problem of scheduling autonomous mobile robots (AMRs) in hospital settings, considering dynamic requests with different priorities. The primary objective is to minimize the daily service cost by dynamically planning routes for the limited number of available AMRs. The total cost consists of AMR's purchase cost, transportation cost, delay penalty cost, and loss of denial of service. To address this problem, we have established a two-stage mathematical programming model. In the first stage, a tabu search algorithm is employed to plan prior routes for all known medical requests. The second stage involves planning for real-time received dynamic requests using the efficient insertion algorithm with decision rules, which enables quick response based on the time window and demand constraints of the dynamic requests. One of the main contributions of this study is to make resource allocation decisions based on the present number of service AMRs for dynamic requests with different priorities. Computational experiments using Lackner instances demonstrate the efficient insertion algorithm with decision rules is very fast and robust in solving the dynamic AMR routing problem with time windows and request priority. Additionally, we provide managerial insights concerning the AMR's safety stock settings, which can aid in decision-making processes.

**1. Introduction**

Given the increasingly severe global aging population and the shortage of medical staff, providing high-quality medical services has become a challenging task. Yen et al. (2018) emphasized that non-nursing activities, such as food tray delivery and equipment retrieval, consume approximately 40% of nurses' time. Furthermore, traditional carts or manual transportation equipment often encounter difficulties when they operate in narrow hospital corridors and crowded areas, impeding the smooth delivery of items. Consequently, hospitals must adopt more efficient and convenient transportation systems to enhance service efficiency and cost savings. Fragapane et al. (2020) compared five innovative applications of autonomous mobile robots (AMRs) with other material handling systems applied in hospitals. The study highlights the benefits of using AMRs in healthcare environments. AMRs can perform various tasks in the medical environment, such as delivering medication, meals, and documents. They can collaborate closely with hospital personnel to enhance the efficiency and quality of medical services by alleviating the workload of medical personnel (Fragapane et al. 2023), thereby providing enhanced convenience for hospital management.

Although AMRs are highly respected for their autonomy and flexibility, they have limitations in terms of loading capacity and battery capacity. In addition, AMRs need to respond quickly to service requests in dynamic medical environments, which poses significant challenges to AMR scheduling and routing strategies. This paper investigates the dynamic AMR routing problem with time windows and requests priority (DAMRRP-TWRP) in the context of AMRs participating in healthcare services.

For AMRs participating in healthcare services, the travel time is uncertain due to the necessity of avoiding obstacles and waiting for elevators while driving. Service time varies due to the different demands of each medical request. Moreover, AMRs need to serve two types of medical requests: static requests and dynamic requests. If all information (such as location, demand, time window, etc.) of requests is known before service starts, the requests are called static requests. Other requests that arrive dynamically during the service process can be called dynamic requests. This study classifies the received dynamic requests into high-priority and low-priority. To prevent excessive costs incurred by introducing too many AMRs into the hospital, we set an upper limit on the total number of AMRs participating in medical services. AMRs have the option to reject servicing low-priority requests.

There are two major difficulties in solving DAMRRP-TWRP. The first challenge is that the dynamic nature changes the information of problems and may affect the corresponding optimal solutions. Therefore, once the change happens, it brings challenges that the initial routing plan might be ineffective and the routes should be replanned with feasible solutions for the new environments. The second difficulty lies in providing high-quality service for dynamic requests with different priorities, given the limited number of available AMRs. Balancing the allocation of resources to cater to these dynamic requests becomes a challenging task. Hence, our study aims to optimize the decision-making process for AMR participation in healthcare services by balancing the total cost of hospitals and customer satisfaction, thus ensuring efficient logistics operations.

In this paper, we established a two-stage mathematical programming model for the DARRP-TWRP. In the first stage, we optimize the allocation of AMRs and their corresponding service routes based on all static requests received at the beginning. Subsequently, in the second stage, the service routes of AMRs are immediately re-optimized based on new requests received. To address large instances and real-time operations, we propose a hybrid meta-heuristic algorithm based on a tabu search (TS) algorithm and efficient insertion algorithm with decision rules (EIADR) and analyze the performance of the algorithm through numerical comparison. The method allows AMRs to wait at the depot and customer locations, as well as strategically accept and reject new customer requests in order to maximize the number of served customers. To the best of our knowledge, this paper represents the first attempt to address the optimization scheduling problem associated with AMR services that involve priority requests and a limited number of AMRs.

The remainder of the paper is structured as follows. Section 2 is devoted to the literature review. In Section 3, we describe the DAMRRP-TWRP problem. In Section 4, we present a mathematical model for the DAMRRP-TWRP. We describe the TS algorithm and EIADR for solving the DAMRRP-TWRP in Section 5. Section 6 presents the computational results. Finally, important conclusions from the results are discussed.

## 2. Literature review

The field of autonomous delivery vehicles and robots has witnessed significant advancements in recent years. AMRs have been increasingly deployed in various intralogistics operations, such as manufacturing, warehousing, cross-docking facilities, terminals, and hospitals. Bhosekar et al. (2021) utilized automated material handling systems to optimize the flow of materials and reduce workforce requirements in medical institutions. The primary research objective centered around two key aspects: path redesign and operational fleet sizing. Cheng et al. (2023a) investigated the optimization problem of scheduling strategies for AMRs at smart hospitals, where the service and travel times of AMRs are stochastic. This paper investigates the dynamic AMR routing problem with time windows and requests priority (DARRP-TWRP), which not only considers the capacity of the AMR, time windows, and demand for medical requests, but also incorporates factors such as traffic conditions, dynamically arriving requests, and request priorities. The DARRP-TWRP is an extension of the dynamic vehicle routing problem. These aspects are briefly discussed in the following.

The dynamic vehicle routing problem (DVRP) is a variant of the vehicle routing problem (VRP) in which the problem inputs are partially known in advance and dynamically change during the operation period (Zhang and Woensel 2023). There are three main aspects of the stochastic elements considered in DVRP including stochastic customer request (Jia et al. 2018; Okulewicz and Mandziuk 2019), service time (Bian and

Liu 2018; Ulmer et al. 2021), and travel time (Aragão et al. 2019; Sabar et al. 2019). The main task of DVRP is to dynamically plan the optimal driving route of these vehicles with the dynamic factors so that these vehicles can meet customers' needs and time windows (Pillac et al. 2013; Abdirad et al. 2021). Sabar et al. (2019) considered the DVRP where the travel time from point to point is influenced by factors like traffic congestion. Bian and Liu (2018) focused on the operational-level stochastic orienteering problem, in which travel time and service time are stochastic. Zhang et al. (2019) optimized the electric vehicle routing problem with stochastic demands by using a traditional recourse strategy and preventive replenishment strategy. Jia et al. (2018) proposed a dynamic logistic dispatching system that allows new orders to be received as the working day progresses. Over the past few decades, most related research has considered the routing problem with one or two dynamic properties. There has been little research on the DVRP with time windows, which simultaneously consider stochastic service time, stochastic travel time, and stochastic demands.

Most research on DVRP uses similar objective functions to that of the classic VRP, such as minimization of total travel time or costs. Hyland and Mahmassani (2018) presented the on-demand shared-use autonomous vehicle mobility service with no shared rides, to minimize travel time and waiting time. Vodopivec and Miller-Hooks (2017) studied DVRP to minimize the costs incurred en route by vehicles. Some scholars have introduced other objectives for DVRP such as expected response time, satisfaction level of the customers, or minimizing delays. Schyns (2015) studies a dynamic capacitated vehicle routing problem with time windows, (partial) split delivery, and heterogeneous fleet. The objective of the proposed approaches is to optimize the responsiveness so that the vehicle can restart its activity as soon as possible. Brinkmann et al. (2019) presented the stochastic-dynamic inventory routing problem for bike-sharing systems and established a reliable system to maximize the users' satisfaction. Ferrucci and Bock (2016) proposed a new proactive real-time routing approach for the dynamic vehicle routing problem with soft time window constraints to minimize total lateness. In the case of hard time windows and a finite number of vehicles, not considering the rejection of customers may render problem instances infeasible (Psaraftis et al. 2016). In fact, it makes sense to reject some customers in some situations. Fikar (2018) studied a dynamic pickup and delivery routing problem that optimally serves all requests within a limited number of vehicles. If a request cannot be delivered due to violating the time window, the order is rejected. Ulmer et al. (2021) considered a stochastic dynamic pickup and delivery problem where the fleet size is limited, time constraints are soft, and drivers are required to accept all orders that are assigned during the order horizon.

To tackle DVRP, researchers and practitioners have often called upon heuristics or hybrid methods. Heuristics are largely employed in solving dynamic problems because they can often produce a good quality solution in modest computational times, but at the expense of optimality (Vidal et al. 2013). The most popular heuristic algorithms are the genetic algorithm (AbdAllah et al. 2017; Achamrah et al. 2022), variable neighborhood search algorithm (Sarasola et al. 2011; Lu et al. 2020), ant colony algorithm (Davis 2017; Xiang et al. 2020), particle swarm optimization algorithms (Okulewicz and Mandziuk 2019), tabu search algorithm (TS) (Hanshar and Ombuki-Berman 2007; Liao and Hu 2011; Gmira et al. 2021), etc. The TS algorithm is an effective method for solving DVRP, which can optimize the routes when an emergency request occurs (Gendreau et al. 2006) and the travel times are uncertain (Cheng et al. 2023b). Compared with other algorithms, this algorithm is the closest to the optimal solution, but the running time is relatively long (Zhang et al. 2022). For this reason, many scholars have made corresponding improvements based on the TS algorithm. Ge et al. (2020) combined Clarke and Wright's saving method with the TS algorithm to solve the electric vehicle routing problem with stochastic demands, which includes two cases where the fixed customers with uncertain demand and random customers. Li and Li (2020) proposed an improved TS algorithm based on a greedy algorithm to solve the vehicle routing problem with stochastic travel time and service time.

When the probability distribution of dynamic input information is unknown, DVRP solutions typically adopt the myopic re-optimization method, which is divided into two methods: periodic re-optimization and instant re-optimization. Periodic re-optimization divides the planning period into distinct time slots or epochs, treating each interval as a separate static routing problem that needs to be solved. For instant re-optimization, once the input data changes, re-optimization is initiated to adjust and optimize the dynamic solution. Hong (2012) proposed a large neighborhood search algorithm for the real-time vehicle routing problem with time windows, in which once a new request arrives, it is immediately considered to be included in the present solution. Ninikas and Minis (2014) proposed an $N$-request re-optimization framework in which the sequential SVRPs are reoptimized when the number of pending dynamic requests reaches a predefined value $N > 1$. Steever et al. (2019) proposed methods embedded in myopic re-optimization to address the dynamic courier routing for a food delivery service.

Compared with previous studies, the current study presents the following theoretical and practical contributions: (1) a two-stage mathematical programming model is established for the dynamic AMR routing problem with time windows and requests priority; (2) a quick response algorithm is developed for the dynamic AMR routing problem; (3) an extensive numerical study that shows valuable insights for safety stock and the service rate for medical requests.

## 3. Problem description

In this section, we first present a formal description of the DAMRRP-TWRP. We then established a two-stage mathematical programming model and provided some decision rules for the DAMRRP-TWRP.

### 3.1. Definition of DAMRRP-TWRP

Let's assume that a hospital deploys $M$ identical AMRs with a capacity of $Q$ for item delivery. We denote the set of routes travelled by the $k$th AMR by $L_k = \{1, 2, \ldots, l_{\max}^k\}$, $k = 1, 2, \ldots, M$. Each AMR operates at a speed of one distance unit per one time unit. The unit travel cost for an AMR is defined as $\xi_1$. All AMRs are fully charged when leaving the electrical charging station at the pharmacy, and the AMRs do not charge during the operation process. While the hospital serves static medical requests (with known information such as location, demand, and time windows), dynamic medical requests with known demand and time windows emerge dynamically over time. The time of receiving a medical request is denoted as $a_i$, with $a_i = 0$ representing the arrival time of static requests.

We define the DAMRRP-TWRP using a complete graph $G = (V, A)$, where $V = \{0, 1, 2, \ldots, n, 0'\}$ represents the set of vertices and $A = \{(i, j) | i, j \in V, i \neq j\}$ represents the set of edges. The vertices 0 and $0'$ represent the depot, while $R = \{1, 2, \ldots, n\} = R_1 \cup R_d$ represents the set of requests, $R_1$ represents the set of static requests and $R_d$ represents the set of dynamic requests. Let $T_{pij}^k$ denote the travel time spent by the $k$th AMR on the route $p \in L_k$ from request $i$ to request $j$. Additionally, we use $S_i$ to denote the service time of request $i$. In our discussion, we assume $T_{pij}^k$ and

$S_i$ are random variables that subject to normal distribution. Each medical request $i \in R$ has its own demand $q_i$ $(q_i \leq Q)$ and time window $[e_i, h_i]$.

All medical requests require the same type of goods, enabling an en-route AMR to be dispatched to serve a newly received request without returning to the depot, provided that it has sufficient capacity. However, it is crucial to acknowledge that accommodating all request demands may not be feasible due to the limited number of available AMRs and the requirement to satisfy both the AMR's capacity and the request time windows. Generally, if a new request consumes too much time or AMR resources and reduces the expected returns, the request is rejected (Zhang and Woensel 2023). In this study, incoming medical requests are classified into high-priority request $i \in R_2$ and low-priority request $i \in R_3$, i.e., $R_d = R_2 \cup R_3$. When the number of AMRs participating in the service has reached the threshold $M$, high-priority requests must be served, even if it leads to a penalty $\xi_2$ for one unit delayed arrival time. For low-priority requests, if the penalty for violating the time window exceeds the loss $\xi_0$ incurred by rejecting the request, the AMR will reject it. Rejected requests are considered lost with corresponding losses $\xi_0$. It is assumed that rejected requests are handled by nurses. Below, we introduce Decision rule 1 on whether to accept services for low-priority requests.

**Decision rule 1** (Acceptance Rule). When the number of AMRs participating in the service has reaches the threshold $M$, the $k$th AMR's route $0 \to i_1 \to \cdots \to i_q \to \cdots \to i_n \to 0'$ accepts and schedules the low-priority request with $i_d \in R_3$ preceding $i_q$ if and only if $\xi_0 > \xi_2 \cdot \left[ \sum_{i=i_d}^{i_n} E\left(A_{pi}^{k\prime} - h_i\right)^+ - \sum_{i=i_q}^{i_n} E\left(A_{pi}^{k} - h_i\right)^+ \right]$ holds, where $A_{pi}^{k}$ is the arrival time at each request before insertion $i_d$, and $A_{pi}^{k\prime}$ is the arrival time at request $i$ after insertion $i_d$.

To better understand the DAMRRP-TWRP, Figure 1 presents an illustrative example of the DAMRRP-TWRP. The example includes a single depot, two AMRs, a set of static requests $\{1,2,\cdots,7\}$, and a set of dynamic requests $\{8,9\}$, where request 8 has high-priority, while request 9 has low-priority. For this example, we assume the present number of AMRs participating in the service is $M=2$. Let $\xi_0=1$, and $\xi_2=1$. Request 8 has a time window of [9:05, 9:15] and can be serviced by the first AMR after completing request 2. On the other hand, request 9 has a time window of [9:10, 9:20]. If the first AMR completes the service of request 8 and goes on to serve requests 9 and 3 sequentially, it would result in an additional delay penalty cost of $\xi_2\left[(9:30-9:20)+(9:50-9:50)\right] > \xi_0$. Therefore, according to Decision rule 1, request 9 is rejected.

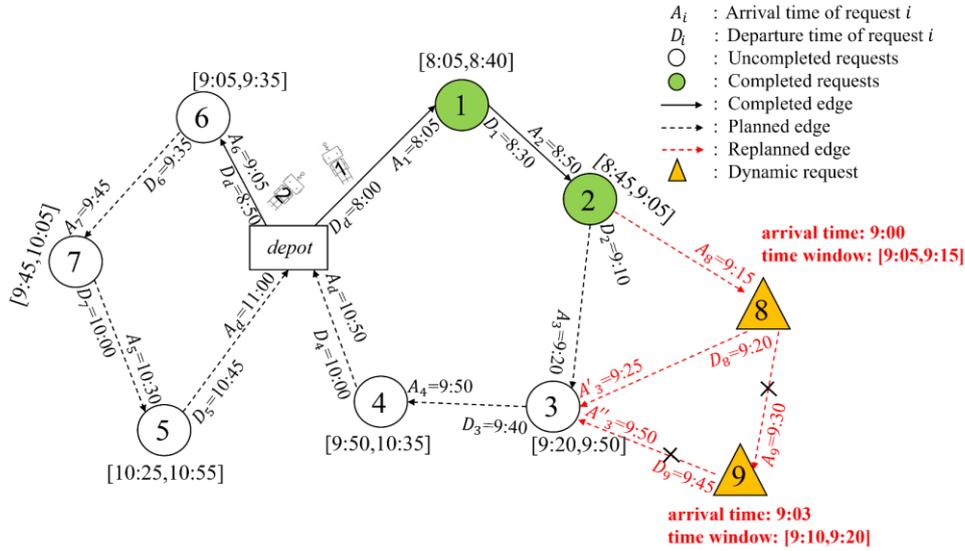

Figure 1. Example of receiving the dynamic request.

In this paper, two strategies are proposed to effectively manage newly arrived requests. Firstly, a safety stock of AMRs denoted by $\psi \in [0,1]$ is introduced. This parameter ensures that the remaining load capacity of each AMR after serving each individual requests remains at least $Q \cdot \psi$. This approach enhances the responsiveness of

AMRs in handling dynamically arriving medical requests. Secondly, if some waiting time is anticipated at the next request location for AMRs, they are instructed to remain stationary at their present position. This strategy allows AMRs to depart from the requesting node at the latest possible moment, thereby facilitating last-minute adjustments to the planned route in response to the influx of new requests.

Notations used in the problem statement are defined as shown in Table 1.

Table 1. Notation for the DAMRRP-TWRP.

| | |
|---|---|
| **Sets** | |
| $R_1 = \{1, 2, \cdots, n_r\}$ | Set of static requests |
| $R_2 = \{1, 2, \cdots, n_h\}$ | Set of high-priority requests |
| $R_3 = \{1, 2, \cdots, n_l\}$ | Set of low-priority requests |
| $L_k$ | Set of the route traveled by the $k$th AMR |
| $X_{pij}^k = \{x_{pij}^k\}$ | Set of $x_{pij}^k$ |
| **Parameters** | |
| $H$ | A big number |
| $M$ | Upper limit of total AMRs |
| $Q$ | Cargo capacity of the AMRs |
| $\xi_0$ | Loss of denial of service |
| $\xi_1$ | Travel cost per unit time |
| $\xi_2$ | Penalty cost per unit time for violating the time window |
| $\xi_3$ | Fixed cost of the AMR |
| $a_i$ | Receival time of the medical request |
| $q_i$ | Demand of request $i$ |
| $S_i$ | Service time of request $i$ |
| $T_{pij}^k$ | Travel time spent by the $k$th AMR on the route $p \in L_k$ from request $i$ to request $j$ |
| $\psi$ | Safety stock |
| $u_{pi}^k$ | Remaining capacity of the $k$th AMR when it arrives at request $i$ on its $p$th route |
| $\gamma_{pi}^k$ | Remaining battery of the $k$th AMR when it arrives at request $i$ on its $p$th route |
| $D_{pi}^k$ | Time when the $k$th AMR departure from request $i$ on the $p$th route |
| $A_{pi}^k$ | Time when the $k$th AMR arrive at the request $i$ on the $p$th route |
| $[e_i, h_i]$ | Time window of request $i$ |
| **Decision variables** | |
| $x_{pij}^k$ | 1 if AMR $k$ serves requests $i$ and $j$ sequentially on its $p$th route; 0 otherwise |
| $y_i$ | 1 if the low-priority dynamic request $i \in R_3$ is served; 0 otherwise |

### 3.2. Model formulation

In this section, we establish a two-stage mathematical programming model for the DAMRRP-TWRP. Firstly, we optimize the service routes of all static requests to establish a set of prior routes. Subsequently, based on these prior routes, optimization adjustments are made for the number of AMRs and their service routes in response to received dynamic requests.

### 3.2.1. Mathematic model for static requests

A static mathematical model is established to address the routing problem of static requests. This problem can be mathematically formulated with parameters $\xi_1$ representing the AMR transportation cost per unit time, $\xi_2$ representing the penalty cost per unit time for violating the time window, and $\xi_3$ representing the fixed cost of the AMR such that $\xi_1 < \xi_2 < \xi_3$. We define the following decision variable: $x_{pij}^k$ equal to 1 if arc $(i, j)$ is traversed by an AMR, 0 otherwise. The model can then be formulated as follows:

$$\min\ \xi_1 \sum_{k=1}^{m} \sum_{p \in L_k} \sum_{\substack{i \in V/(R_2 \cup R_3) \\ j \in V/(R_2 \cup R_3) \\ i \neq j}} E\left(T_{pij}^k\right) \cdot x_{pij}^k + \xi_2 \sum_{k=1}^{m} \sum_{p \in L_k} \sum_{i \in R_1} E\left(A_{pi}^k - h_i\right)^+ + \xi_3 m \tag{1}$$

$$\text{s.t.}\ \sum_{k=1}^{m} \sum_{p \in L_k} \sum_{j \in V, i \neq j} x_{pij}^k = 1,\ i \in R_1 \tag{2}$$

$$\sum_{j \in V/\{0\}, j \neq i} x_{pij}^k = \sum_{j \in V/\{0'\}, j \neq i} x_{pji}^k,\ i \in V/R_d, k = 1, 2, \ldots, m, p \in L_k \tag{3}$$

$$u_{p0}^k = Q,\ k = 1, 2, \ldots, m, p \in L_k \tag{4}$$

$$Q\psi \leq u_{pi}^k - q_i \leq Q,\ i \in V/R_d, k = 1, 2, \cdots m, p \in L_k \tag{5}$$

$$q_j \leq u_{pj}^k \leq u_{pi}^k - x_{pij}^k \cdot q_i + H(1 - x_{pij}^k), i \in \{0\} \cup R_1, j \in R_1 \cup \{0'\}, i \neq j, \quad (6)$$
$$k = 1, 2, \cdots, m, p \in L_k$$

$$A_{p,0'}^k = A_{p+1,0}^k, k = 1, 2, \cdots, m, p \in L_k \quad (7)$$

$$D_{pi}^k = max\left\{A_{pi}^k + S_i, e_j - T_{pij}^k\right\}, i, j \in V/R_d, k = 1, 2, \cdots m, p \in L_k, x_{pij}^k = 1 \quad (8)$$

$$A_{pj}^k = D_{pi}^k + T_{pij}^k, i, j \in V/R_d, k = 1, 2, \cdots m, p \in L_k, x_{pij}^k = 1 \quad (9)$$

$$A_{pi}^k + S_i + T_{pij}^k \cdot x_{pij}^k - H \cdot \left(1 - x_{pij}^k\right) \leq A_{pj}^k, i \in \{0\} \cup R_1, j \in R_1 \cup \{0'\}, i \neq j, \quad (10)$$
$$k = 1, 2, \cdots, m, p \in L_k$$

$$x_{pij}^k \in \{0,1\}, i \in \{0\} \cup R_1, j \in R_1 \cup \{0'\}, i \neq j, k = 1, 2, \ldots, m, p \in L_k \quad (11)$$

The objective function (1) aims to minimize the total cost. The first term corresponds to the expected transportation cost, the second term relates to the loss for expected delay time, and the third term is the fixed cost of AMRs. Constraint (2) guarantees each request to be visited exactly once by an AMR. Constraint (3) indicates that the number of input arcs at each point is equal to the number of output arcs. The constraints (4)-(6) represent the remaining loading capacity constraints of the AMR on any route. The constraints (7)-(10) define the departure time and arrival time of AMRs. Equation (8) ensures that when an AMR arrives at request $j$, it initiates service directly without waiting, thus achieving a "last-minute" route change. The constraint (11) denotes a binary decision variable $x_{pij}^k$. Through this formulation, a set of priori routes for AMRs can be derived. These planned routes are not changed or updated after realization unless dynamic requests arrive.

We take the approach of Ehmke et al. (2015) to approximate the arrival time by a normal distribution, i.e., $A_{pi}^k \sim N\left(\mu\left(A_{pi}^k\right), \sigma^2\left(A_{pi}^k\right)\right)$, where $A_{pj}^k = D_{pi}^k + T_{pij}^k$ if $x_{pij}^k = 1, \forall i, j \in V$. The expected penalty for delayed service is

$$\xi_2 \cdot E\left(A_{pi}^k - h_i\right)^+ = \xi_2 \cdot E\left(\max\left\{0, A_{pi}^k - h_i\right\}\right) = \xi_2 \cdot E\left(\max\left\{h_i, A_{pi}^k\right\} - h_i\right), \quad (12)$$

The departure time $D_{pi}^k$ can be expressed as $D_{pi}^k = \max\{A_{pi}^k + S_i + T_{p,i,i+1}^k, e_{i+1}\} - T_{p,i,i+1}^k$.

Since $A_{p,i+1}^k = D_{p,i}^k + T_{p,i,i+1}^k$, we have $A_{p,i+1}^k = \max\{A_{pi}^k + S_i + T_{p,i,i+1}^k, e_{i+1}\}$. Refer to Nadarajah and Kotz (2008), we have

$$\mu(A_{p,i+1}^k) = \mu_1 \Phi\left(\frac{\mu_1 - e_{i+1}}{\sigma_1}\right) + e_{i+1} \Phi\left(\frac{e_{i+1} - \mu_1}{\sigma_1}\right) + \sigma_1 \cdot \phi\left(\frac{\mu_1 - e_{i+1}}{\sigma_1}\right),$$

$$\sigma^2(A_{p,i+1}^k) = (\sigma_1^2 + \mu_1^2) \cdot \Phi\left(\frac{\mu_1 - e_{i+1}}{\sigma_1}\right) + e_{i+1}^2 \cdot \Phi\left(\frac{e_{i+1} - \mu_1}{\sigma_1}\right)$$

$$+ (\mu_1 + e_{i+1})\sigma_1 \cdot \phi\left(\frac{\mu_1 - e_{i+1}}{\sigma_1}\right) - \mu^2(A_{p,i+1}^k),$$

where $\mu_1 = \mu(A_{pi}^k) + \mu(S_i) + \mu(T_{p,i,i+1}^k)$, $\sigma_1 = \sigma(A_{pi}^k + S_i + T_{p,i,i+1}^k)$, $\Phi(\cdot)$ is the cumulative distribution function of the standard normal distribution, and $\phi(\cdot)$ is the probability density function of the standard normal distribution. Therefore, according to $\mu(A_{pi}^k)$ and $\sigma^2(A_{pi}^k)$, equation (12) can be expressed as

$$\xi_2 \cdot E(A_{pi}^k - h_i)^+ = \xi_2 \cdot E(\max\{h_i, A_{pi}^k\} - h_i)$$

$$= \xi_2 \cdot \left( h_i \Phi\left(\frac{h_i - \mu(A_{pi}^k)}{\sigma(A_{pi}^k)}\right) + \mu(A_{pi}^k) \Phi\left(\frac{\mu(A_{pi}^k) - h_i}{\sigma(A_{pi}^k)}\right) + \sigma(A_{pi}^k) \phi\left(\frac{h_i - \mu(A_{pi}^k)}{\sigma(A_{pi}^k)}\right) - h_i \right).$$

*3.2.2. Mathematic model for dynamic requests*

The objective of the dynamic programming stage is to efficiently modify a set of priori routes obtained by the model in Section 3.2.1, based on newly received dynamic requests. The goal is to minimize the cost associated with the adjusted routes. It should be noted that due to the upper limit of AMRs participating in medical service, a feasible solution may not always be attainable. However, in this paper, to obtain an effective solution for the hospital, refusing service of low-priority dynamic requests is allowed, which causes additional costs. Hence minimizing the loss of denial of service is incorporated into the objective function, while the AMR needs to meet loading capacity and time window constraints. Let $C_i$ denote the additional cost associated with inserting dynamic request $i \in R_2 \cup R_3$, which consists of the transportation cost, delay penalty cost, and the fixed cost of AMRs. Additionally, we define parameter $y_i$ equal to 1 if the low-priority dynamic request $i \in R_3$ is served, 0 otherwise.

The model can then be formulated as follows:

$$\min \sum_{i \in R_2} C_i + \sum_{i \in R_3} \left( \xi_0 (1 - y_i) + C_i y_i \right) \tag{13}$$

$$\text{s.t.} \sum_{k=1}^{m} \sum_{p \in L_k} \sum_{j \in V, i \neq j} x_{pij}^k = 1, \ i \in R_2 \tag{14}$$

$$\sum_{j \in V/\{0\}, j \neq i} x_{pij}^k = \sum_{j \in V/\{0'\}, j \neq i} x_{pji}^k, \ i \in V, k = 1, 2, \ldots, m', p \in L_k \tag{15}$$

$$u_{p,0}^k = Q, \ k = 1, 2, \ldots, m', p \in L_k \tag{16}$$

$$q_j \leq u_{pj}^k \leq u_{pi}^k - x_{pij}^k \cdot q_i + H(1 - x_{pij}^k), \ i \in V/\{0'\}, j \in V/\{0\}, i \neq j,$$
$$k = 1, 2, \cdots, m', p \in L_k \tag{17}$$

$$A_{p,0'}^k = A_{p+1,0}^k, \ k = 1, 2, \cdots, m', p \in L_k \tag{18}$$

$$D_{pi}^k = \max \left\{ A_{pi}^k + S_i, e_j - T_{pij}^k \right\}, \ i, j \in V, k = 1, 2, \cdots m', p \in L_k, x_{pij}^k = 1 \tag{19}$$

$$A_{pj}^k = D_{pi}^k + T_{pij}^k, \ i, j \in V, k = 1, 2, \cdots m', p \in L_k, x_{pij}^k = 1 \tag{20}$$

$$A_{pi}^k + S_i + T_{pij}^k \cdot x_{pij}^k - H \cdot (1 - x_{pij}^k) \leq A_{pj}^k, \ i \in V/\{0'\}, j \in V/\{0\}, i \neq j,$$
$$k = 1, 2, \cdots, m', p \in L_k \tag{21}$$

$$x_{pij}^k \in \{0,1\}, \ i \in V/\{0'\}, j \in V/\{0\}, i \neq j, k = 1, 2, \ldots, m', p \in L_k \tag{22}$$

$$y_i \in \{0,1\}, i \in R_3 \tag{23}$$

$$m \leq m' \leq M \tag{24}$$

The objective function (13) aims to minimize the total extra cost. The first term corresponds to the additional cost of inserting high-priority requests, including the transportation cost, the delay penalty cost, and the fixed cost of AMRs. The second term relates to the service cost of low-priority requests, including the loss of denial of service and additional costs incurred by the service. Constraint (14) guarantees each high-priority request to be visited exactly once by an AMR. Constraint (15) indicates that the number of input arcs at each point is equal to the number of output arcs. The constraints (16)-(17) represent the remaining loading capacity constraints of the AMR on any route. The constraints (18)-(21) define the departure time and arrival time of AMRs. The constraints (22)-(23) denote binary decision variable $x_{pij}^k$ and $y_i$. The constraint (24) represents the number of AMRs participating in service.

**Definition 1.** For any route $0 \to \cdots \to i \to j \to \cdots 0'$, if the total actual demand on the route $0 \to \cdots \to i \to j$ exceeds $Q$ and the total actual demand on the route $0 \to \cdots \to i$ does not exceed $Q$, we call it a **failure route** due to customer $j$.

**Definition 2.** A route $0 \to \cdots \to i \to j \to \cdots 0'$ is called an **infeasible route** if the time window constraints not satisfied.

For a newly arrived dynamic request $i \in R_2 \cup R_3$, there may be multiple service schemes $X_{pij}^k \neq \varnothing$. We propose an assignment rule (Decision rule 2) for dynamic requests, which aims to schedule new requests in a manner that minimizes the additional cost of adjusting the route. The specific criteria are defined as follows:

**Decision rule 2** (Assignment Rule). Without considering route failure (Definition 1), the optimal assignment rule for request $i \in R_2 \cup R_3$ would be either

$$x_{pij}^{k*} = \arg\min_{x_{pij}^k \in X_{pij}^k} C_i\left(x_{pij}^k, \Delta Travel, \Delta Time, \Delta AMR\right), \forall i \in R_2$$

or

$$x_{pij}^{k*} = \arg\min_{x_{pij}^k \in X_{pij}^k} \left(\xi_0(1 - y_i) + C_i\left(x_{pij}^k, \Delta Travel, \Delta Time, \Delta AMR\right) \cdot y_i\right), \forall i \in R_3.$$

Here, $C_i\left(x_{pij}^k, \Delta Travel, \Delta Time, \Delta AMR\right)$ represents the additional cost associated with assigning dynamic request $i \in R_2 \cup R_3$ before the request $j$ on the $p$th path of the $k$th

AMR. This additional cost includes the extra transportation cost $\Delta Travel$, the extra delay penalty cost $\Delta Time$, and the extra fixed cost of AMRs $\Delta AMR$.

## 4. Algorithm for the DAMRRP-TWRP

The DARRP-TWRP is an extension of the dynamic vehicle routing problem. Consequently, the DAMRRP-TWRP falls into the category of NP-hard problems as well. When the problem size or complexity exceeds the capabilities of CPLEX, a common approach is to employ heuristics or metaheuristics, such as those discussed in Devapriya et al. (2017), to find the best possible solution. This section describes the metaheuristic algorithm proposed to solve the DAMRRP-TWRP and presents its framework. The algorithm is divided into two stages. In the first stage, a Tabu Search framework is used to optimize the routes and number of AMRs servicing all static requests (see Section 4.1). The newly dynamic request information in the hospital is unknown until the dynamic request arrives. To address this, an instant re-optimization strategy is considered, which involves re-optimizing the prior routes based on the dynamic requests that arise. In the second stage, we propose a quick response algorithm named efficient enhanced insertion algorithm with decision rules (EIADR) to optimize the number and routes of AMRs after the arrival of dynamic service requests (see Section 4.2).

### *4.1. TS framework to optimize static requests*

To solve the mathematical model in Section 3.2.1, we adopt Tabu Search (TS) to plan prior routes for all static requests. The Tabu Search (TS) was proposed by Glover (1986) which explores part of the solution space by moving to the best neighbor of the present solution, even when this movement deteriorates the objective function. This approach is a powerful and effective heuristic algorithm for solving vehicle routing problems. An overview of the main components of the algorithm is presented in the following sections.

*4.1.1. Initial solution*

For heuristic algorithms, the quality of the initial solution directly affects whether the algorithm can efficiently find the optimal solution to the problem. Therefore, it is necessary to choose a reliable method. Here, we consider using an improved greedy algorithm to generate the initial solution of TS, as detailed in Algorithm A.1 (Appendix A).

*4.1.2. Neighborhood structure*

To obtain improved solutions within a limited time, we use three enhanced operators proposed by Cheng et al. (2023b), namely swap*, 2-opt*, and relocation*. The selection of the operators follows a roulette wheel procedure where the probability of selecting a certain operator depends on its weight. These weights are updated according to a reward scheme.

The solutions generated by the neighborhood operators might be infeasible because they do not consider load and time window constraints. Hence, we introduce two repair operators.

**Repair operator 1:** If a service route failure occurs due to overloading, a depot is inserted before the point where the remaining load does not meet the downward load demands.

**Repair operator 2:** If a route is infeasible due to violating time windows, an AMR is inserted before the first node where the time window is not satisfied.

*4.1.3. Tabu list and tabu tenure*

The tabu mechanism in the TS algorithm can explore a wide range of solutions and avoid getting stuck in local optima in the search process. Here we use three tabu lists $B_i, \forall i = 1, 2, 3$ to store the solutions produced by applying each of the three neighborhood operators to each pair of requests, where the tabu list $B_i$ is an $n \times n$ symmetric matrix, and $n$ is the total number of requests. For the operator $i$, the element $B_i(j_1, j_2)$ in the tabu list $B_i$ represents the tabu state of the neighborhood action $b_i(j_1, j_2)$, where $j_1, j_2 \in R$ are the two requests for swap, inversion or relocation. As done in Gendreau et al. (1996), we consider a random tabu tenure, which is an integer uniformly generated from $[L_1, L_2]$.

*4.1.4. Aspiration criterion*

To prevent the algorithm from prematurely terminating in a suboptimal solution, we consider using the aspiration criterion. For a neighborhood solution that is currently tabu but better than the present optimal solution, we revoke the tabu operation and take the solution as the present optimal solution.

*4.1.5. Algorithm stop condition*

In this paper, for the stopping condition of the TS algorithm, we consider that the algorithm stops after a predetermined maximum number of iterations $N$ and outputs the final solution.

The TS algorithm proposed in this paper to solve the DAMRRP-TWRP is summarized in Algorithm A.2 (Appendix A).

*4.2. Enhanced insertion algorithm with decision rules for dynamic requests*

To make a quick response for the dynamic requests and get an approximate optimal solution to the mathematical model in Section 3.2.2, this subsection proposes an enhanced insertion algorithm with decision rules (EIADR) to handle dynamically arrived requests. This algorithm is based on a myopic re-optimization of the AMR routes whenever one or multiple requests arise. The myopic re-optimization strategy transforms dynamic problems into multiple discrete static problems for solution, which can be achieved through two approaches: periodic re-optimization and instant re-optimization. Periodic re-optimization divides the planning period into distinct time slots or epochs, treating each interval as a separate static routing problem that needs to be solved. Decision epochs occur at fixed intervals, but new requests must wait until the end of the current time interval before being scheduled or rejected. If the decision period is too long, it may result in the time window of dynamic requests received during that period not being satisfied. Therefore, our method adopts an instant re-optimization strategy, where the routes are immediately re-planned upon the arrival of new requests. This approach prevents time window violations and ensures the quality of the solution to the best extent possible. In this instant re-optimization, the requests already visited are removed from the present solution. This helps to adjust the routes effectively and incorporate the newly arrived requests.

Dynamic request $i_d \in R_2 \cup R_3$ has multiple insertion schemes on the route of each AMR participating in service, including some infeasible schemes that violate constraints. To quickly and effectively identify feasible insertion positions for dynamic requests on the route set, we introduce Property 1 to quickly eliminate infeasible solutions.

**Property 1.** Let the receiving time of the dynamic request $i_d$ be $a_{i_d}$, its time window be $\left[e_{i_d}, h_{i_d}\right]$. The $k$ th AMR service route $r^k$ is $0 \to i_1^k \to i_2^k \to \cdots \to i_{n-1}^k \to i_n^k \to 0'$. If there is a request $i_j^k$ on route $r^k$ satisfying $P\left\{A_{pi_j^k}^k + S_{i_j^k}^k + T_{i_j^k, i_d}^k < h_{i_d}, a_{i_d} < D_{pi_j^k}^k\right\} > 1 - \varepsilon$, request $i_d$ can be served after $i_j^k$ without breaking its time window with the probability $1 - \varepsilon$.

Based on Property 1, Algorithm 1 presents pseudocode for searching for insertion schemes for dynamic requests in the service routes of $m$ AMRs. Firstly, initialize vector $IP = \varnothing$ to denote the insertion locations of dynamic requests in the service routes of $m$ AMRs, upon the arrival of the dynamic request (line 1). Next, in lines 2-9, based on the proposed Property 1, determine the insertion locations for dynamic request $i_d$ on all the available service routes and add the corresponding scheme to the $ip_k$ (line 5). Finally, all feasible insertion schemes are stored within $IP$ (line 8).

**Algorithm 1** Insertion scheme

1: Initialize insertion location $IP = \varnothing$ on the route of each AMR for the dynamic request $i_d$, and $ip_k = \varnothing$, $\forall ip_k \in IP$. When the dynamic request $i_d$ arrives, the route set of $m$ AMRs is $\mathcal{R} = \{r^1, \ldots, r^k, \ldots, r^m\}$.

2: **For** each AMR's service route $r^k$ is $i_1^k \cdots \to i_j^k \to i_{j+1}^k \to \cdots \to i_n^k \to 0'$ **do**

3:    **For** each request on route $r^k$ **do**

4:       **If** $P\left\{A_{pi_j^k}^k + S_{i_j^k} + T_{i_j^k, i_d}^k < h_{i_d}, a_{i_d} < D_{pi_j^k}^k\right\} > 1 - \varepsilon$ **then** // **Property 1**

5:         $ip_k = ip_k \cup \{i_j^k\}$ // update $ip_k$

6:       **End**

7:    **End**

8:    $IP = IP \cup ip_k$ // Update $IP$

9: **End**

Algorithm 2 presents a comprehensive outline of the EIADR method. For all dynamic requests received by the system, Algorithm 1 is employed to identify all insertion schemes for the dynamic request $i_d$ on the set of routes $\mathcal{R}$ (line 2). If the insertion of the request $i_d$ leads to route failures due to overloading, Repair operator 1 is applied to obtain a set of routes that satisfy the capacity constraints (lines 3-5). If the number of AMRs participating in the service has not reached the upper limit $M$, and the insertion of the request $i_d$ results in infeasible routes due to time window violations, Repair operator 2 is employed to add AMR (lines 7-9). However, if the present number of AMRs has reached the upper limit of $M$, the decision to serve the request $i_d$ depends on its priority. For low-priority requests, apply Decision rule 1 to determine whether to serve the request (line 13). Conversely, if the request $i_d$ has a high priority, Decision rule 2 is utilized to identify the optimal insertion scheme $x^{k*}_{p,i_d,j}$ (line 15). Finally, the route set $\mathcal{R}$ is updated and output (line 20).

---

**Algorithm 2** EIADR for the dynamic request

---

1: **For** a dynamic request $i_d$ with arrival time $a_{i_d}$ and the time window $\left[e_{i_d}, h_{i_d}\right]$ **do**

2: According to **Algorithm 1**, find all insertion schemes on the route set $\mathcal{R}$.

3: **If** the dynamic request $i_d$ is inserted into the route of the $k$th AMR, the route is failure due to overloading **then**

4:   Correct the route to a feasible route according to the **Repair operator 1**.

5: **End If**

6: **If** the present number of AMR used is less than $M$

7:   **If** the dynamic request $i_d$ is inserted into the route of the $k$th AMR, but the route is infeasible due to violating time windows **then**

8:     Correct the route to a feasible route according to the **Repair operator 2**

9:   **End If**

10: **End If**

11: **If** the present number of AMR used is $M$

12:   **If** the dynamic request $i_d$ is low priority **then**

13:     Determine whether to serve the request $i_d$ based on **Decision rule 1.**

14:   **Else**

15:     Find the best assignment scheme $x^{k*}_{p,i_d,j}$ based on **Decision rule 2.**

16:   **End If**

17: **End If**

18: Update present solution $\mathcal{R}$.

19: **End For**

20: Output the final solution $\mathcal{R}$.

---

Algorithm 3 provides a frame of the metaheuristic algorithm proposed to solve the DAMRRP-TWRP. Use the TS algorithm (see Algorithm A.2) to optimize the routes of all static requests in the AMRs service and obtain the set of priori routes $\mathcal{R}$ (line 1). Lines 2-5 plan all received dynamic requests, whereas line 3 utilizes the EIADR method to plan each received dynamic request based on $\mathcal{R}$ (see Algorithm 2). For situations with low timeliness requirements, the TS algorithm can be chosen to further optimize the solution obtained from EIADR (line 4). We call EIADR+TS method. Whether to perform t EIADR+TS method can refer to Section 5.3.1 of the experiment. This dynamic process will repeat until there are no new requests or services.

---
**Algorithm 3** Dynamic Optimization Algorithm
---
1: Execute TS algorithm on all static requests and create a solution $\mathcal{R}$ as a set of priori routes for AMRs
2: **While** there are dynamic requests
3:   Add new arrival requests to $\mathcal{R}$ by executing EIADR method // **Algorithm 2**
4:   Execute the TS algorithm to further optimize the solution obtained from EIADR //
       EIADR+TS method may or may not be executed.
5: **End While**
---

## 5. Computational experiments

This section describes the computational experiments. The test work is completed on a computer with an Inter(R) Core (TM) i9-9900K CPU @ 3.60GHz, 16.00GB, and an operating system of Windows 10. The algorithm was coded in MATLAB 2018b. Section 5.1 provides a detailed description of the benchmark data. The parameter setting is discussed in Section 5.2. The results of detailed and extensive computational experiments are then presented in Section 5.3.

### *5.1. Benchmark data*

Since the DAMRRP-TWRP is a new problem, there are no instances for it. To evaluate the effectiveness of the proposed model and algorithm, we use the adjusted Lackner's benchmark. The coordinates, demand, and time windows of vertices are the same as in the Lackner instances. The mean service time $E[S_i]$ of each vertex equal to the vertex service time in the Lackner instances. Additionally, the mean values of travel times between nodes are equal to the corresponding Euclidean distances. There are 56 benchmark instances divided into six types: R1, R2, C1, C2, RC1, and RC2, with each type consisting of 8-12 instances. These instances vary across different types based on the distribution of nodes, service time for each node, and width of time windows. Notably, each instance includes five distinct degrees of dynamism, namely 10%, 30%, 50%, 70%, and 90% levels. The dynamic degree of the instance is characterized by:

$$\Delta = \frac{n_d}{n_r + n_d} * 100\%$$

where $n_d = n_h + n_l$ denotes the number of dynamic requests including the number of high priority requests $n_h$ and the number of low priority requests $n_l$, and $n_r$ denotes the total number of static requests. To incorporate various degrees of dynamism, the initial set of 56 base instances is expanded into a total of 280 instances. Throughout the subsequent experiments, each instance is labeled as "data_d%", where "data" represents the base instance data and "d%" denotes the corresponding dynamic degree. For instance, "C101_10%" signifies that it originates from the C101 Lackner instance with a dynamic degree of 10%, meaning there are 10% dynamic requests within the instance.

### 5.2. Parameter setting

The parameters $\xi_0$, $\xi_1$, $\xi_2$, $\xi_3$ are 1000, 1, 100, and 3000 respectively in the objective function. The variances of service time and travel time are 10. For the tabu tenure in the TS algorithm, if the tabu tenure is too long, it is easy to miss the optimal solution, and if it is too short, it is easy to fall into the local optima. Therefore, the choice of the tabu tenure is significant for this algorithm. We assume that the tabu tenure is uniformly generated from $[L_1, L_2] = [40, 50]$. In addition, for safety stock $\psi \in [0,1]$, the experiment is conducted with a value of 0.2. For the total number of AMRs participating in medical services $M$, we set the ratio of AMRs based on the ratio of static requests to total requests, i.e. $\frac{m}{M} = \frac{n_r}{n_r + n_d} \Rightarrow M = \left\lceil \frac{m}{1-\Delta} \right\rceil$, the parameter $m$ is the number of AMR participating in static requests and $\Delta$ is the dynamic degree. Table 2 shows the parameter settings in the experiments.

Table 2. Parameter settings.

| Parameter | Value |
| --- | --- |
| Loss of denial of service for low priority request $\xi_0$ | 1000\$ |
| the unit travel cost of the AMR $\xi_1$ | 1\$/m |
| the unit loss of violating the time window $\xi_2$ | 100\$ |
| the daily fixed cost of each AMR $\xi_3$ | 3000\$ |
| the variance service time $D[S_i]$ | 10 |
| the variance travel time $D[T_{pij}^k]$ | 10 |
| safety stock $\psi$ | 0.2 |
| tabu tenure $[L_1, L_2]$ | $[40, 50]$ |
| iterations $N$ | 500 |

## 5.3. Computational results

This section reports on the numerical results for the adjusted Lackner's benchmark. The effectiveness of the model and relevant algorithms developed in this study is evaluated in Section 5.3.1, some computational results for instances are presented in Section 5.3.2, and several management insights through sensitivity are provided in Section 5.3.3.

*5.3.1. Performance evaluation of the EIADR method*

Because we have changed some settings of the Lackner instances, there is no available proven optimal solution for our problem. No comparative data and no competing heuristics exist for our problem, and it is also hard to compute optimal solutions or tight lower bounds. Therefore, in this sub-section, the EIADR method and its combination with TS are directly used to compare and solve the five dynamic degrees separately. Our computational results are presented in Table 3, Table 4, Table 5, and Table 6. Here we only present the results of all instances with a dynamic degree of 10%. The other four dynamic degrees are presented separately through the average values of each instance. Table 3 and Table 4 report the computational results for the modified type 1 (R1, C1, RC1) and type 2 (R2, C2, RC2) instances with $\Delta = 10\%$, respectively.

For description convenience, we call the column (1) C1, column (2) C2 and so on. C2, C3, C4, and C5 are the computational results for proposed EIADR method which are the average objective value of the solutions found in 10 runs (Obj_A), the average computational time in seconds required to find the solutions in 10 runs (CPU_A), the best objective value of the solutions found in 10 runs (Obj_B), and the best computational time in seconds required to find the solutions in 10 runs (CPU_B), respectively. Similarly, C6, C7, C8, and C9 are the computational results for the EIADR+TS algorithm, which means that after processing all dynamic requests using the EIADR method, the TS algorithm is used to further optimize the solution obtained from EIADR. C10 is the relative error between C4 and C8 (i.e., C10=(C4-C8)/C8*100%), and C11 is the one between C2 and C6 (i.e., C11=(C2-C6)/C6*100%). It represents the quality of the solution generated by directly planning dynamic requests using the EIADR method proposed in this paper.

Table 3. Comparison of EIADR and EIADR+TS method for type 1 instances with

$\Delta = 10\%$.

| Instance (1) | EIADR | | | | EIADR +TS | | | | BRE (10) | ARE (11) |
|---|---|---|---|---|---|---|---|---|---|---|
| | Obj_A (2) | CPU_A (3) | Obj_B (4) | CPU _B (5) | Obj_A (6) | CPU _A (7) | Obj_B (8) | CPU _B (9) | | |
| C101 | 272126.4 | 0.5 | 272126.4 | 0.234 | 272126.4 | 164.98 | 272126.4 | 162.37 | 0 | 0 |
| C102 | 138597.7 | 0.368 | 54418.9 | 0.25 | 137224.0 | 167.61 | 54418.9 | 163.42 | 0 | 1.001 |
| C103 | 41162.1 | 0.45 | 38703.7 | 0.375 | 41162.1 | 164.46 | 38703.7 | 163.26 | 0 | 0 |
| C104 | 36918.7 | 0.409 | 34744.7 | 0.312 | 33857.1 | 145.31 | 29760.7 | 137.89 | 16.74 | 9.042 |
| C105 | 59427.6 | 0.237 | 59427.6 | 0.187 | 59427.6 | 163.72 | 59427.6 | 159.04 | 0 | 0 |
| C106 | 96403.3 | 0.259 | 69619.0 | 0.203 | 96403.3 | 163.92 | 69619.0 | 162.11 | 0 | 0 |
| C107 | 49111.6 | 0.2 | 46180.6 | 0.171 | 49111.6 | 161.74 | 46180.6 | 159.59 | 0 | 0 |
| C108 | 52725.4 | 0.321 | 41752.2 | 0.265 | 52725.4 | 163.25 | 41752.2 | 160.65 | 0 | 0 |
| C109 | 38407.5 | 0.271 | 37305.5 | 0.203 | 38407.5 | 164.04 | 37305.5 | 162.93 | 0 | 0 |
| **Average** | **87208.92** | **0.335** | **72697.62** | **0.244** | **86716.1** | **162.11** | **72143.8** | **159.02** | **0.767** | **0.568** |
| R101 | 275686.5 | 0.303 | 275686.5 | 0.203 | 275686.5 | 164.40 | 275686.5 | 163.34 | 0 | 0 |
| R102 | 214168 | 0.193 | 214168 | 0.140 | 214168 | 163.57 | 214168 | 162.57 | 0 | 0 |
| R103 | 163645.8 | 0.262 | 163645.8 | 0.125 | 163645.8 | 162.50 | 163645.8 | 160.25 | 0 | 0 |
| R104 | 116851.9 | 0.425 | 75083.41 | 0.265 | 116851.9 | 160.53 | 75083.41 | 159.20 | 0 | 0 |
| R105 | 80258.37 | 0.196 | 80258.37 | 0.156 | 80258.37 | 166.10 | 80258.37 | 163.25 | 0 | 0 |
| R106 | 77482.67 | 0.325 | 77482.67 | 0.203 | 77482.67 | 164.9 | 77482.67 | 162.79 | 0 | 0 |
| R107 | 79466.1 | 0.284 | 79466.1 | 0.25 | 79466.1 | 164.75 | 79466.1 | 163.45 | 0 | 0 |
| R108 | 74247.42 | 0.362 | 74247.42 | 0.343 | 74247.42 | 164.42 | 74247.42 | 161.92 | 0 | 0 |
| R109 | 79177.02 | 0.25 | 79177.02 | 0.171 | 79177.02 | 160.92 | 79177.02 | 159.15 | 0 | 0 |
| R110 | 86211.42 | 0.246 | 86211.42 | 0.218 | 86211.42 | 163.45 | 86211.42 | 162.15 | 0 | 0 |
| R111 | 79521.58 | 0.296 | 79521.58 | 0.25 | 79521.58 | 164.54 | 79521.58 | 163.73 | 0 | 0 |
| R112 | 82105.1 | 0.268 | 82105.1 | 0.234 | 82105.1 | 161.36 | 82105.1 | 158.67 | 0 | 0 |
| **Average** | **117401.8** | **0.284** | **113921.1** | **0.213** | **117401.8** | **163.45** | **113921.1** | **161.70** | **0** | **0** |
| RC101 | 72927.09 | 0.15 | 72927.09 | 0.125 | 72927.09 | 143.03 | 72927.09 | 142.59 | 0 | 0 |
| RC102 | 84614.93 | 0.218 | 84614.93 | 0.171 | 84614.93 | 147.11 | 84614.93 | 145.06 | 0 | 0 |
| RC103 | 83003.09 | 0.275 | 83003.09 | 0.218 | 83003.09 | 148.13 | 83003.09 | 146.67 | 0 | 0 |
| RC104 | 63348.10 | 0.271 | 52067.48 | 0.234 | 63348.10 | 137.64 | 52067.48 | 133.67 | 0 | 0 |
| RC105 | 113234 | 0.140 | 113234 | 0.109 | 113234 | 137.56 | 113234 | 137.18 | 0 | 0 |
| RC106 | 75784.05 | 0.215 | 75784.05 | 0.187 | 75784.05 | 137.6 | 75784.05 | 136.96 | 0 | 0 |
| RC107 | 81694.6 | 0.228 | 81694.6 | 0.203 | 81694.6 | 134.16 | 81694.6 | 133.75 | 0 | 0 |
| RC108 | 86810.96 | 0.265 | 86810.96 | 0.234 | 86810.96 | 137.34 | 86810.96 | 136.32 | 0 | 0 |
| **Average** | **82677.10** | **0.220** | **81267.02** | **0.185** | **82677.10** | **140.32** | **81267.02** | **139.02** | **0** | **0** |
| **Avg All** | **95762.61** | **0.279** | **89295.25** | **0.214** | **95598.33** | **155.29** | **89110.64** | **153.25** | **0.256** | **0.189** |

Table 4. Comparison of EIADR and EIADR+TS method for type 2 instances with $\Delta = 10\%$.

| Instance (1) | EIADR | | | | EIADR +TS | | | | BRE (10) | ARE (11) |
|---|---|---|---|---|---|---|---|---|---|---|
| | Obj_A (2) | CPU_A (3) | Obj_B (4) | CPU _B (5) | Obj_A (6) | CPU _A (7) | Obj_B (8) | CPU _B (9) | | |
| C201 | 27596.62 | 0.175 | 23611.63 | 0.140 | 27596.62 | 143.19 | 23611.63 | 141.31 | 0 | 0 |
| C202 | 17493.19 | 0.165 | 17493.19 | 0.109 | 17493.19 | 162.16 | 17493.19 | 159.42 | 0 | 0 |
| C203 | 18897.28 | 0.403 | 14092.15 | 0.328 | 18897.28 | 161.93 | 14092.15 | 160.14 | 0 | 0 |
| C204 | 14997.14 | 0.725 | 13787.79 | 0.640 | 14997.14 | 160.28 | 13787.79 | 156.51 | 0 | 0 |
| C205 | 15647.32 | 0.178 | 13326.43 | 0.125 | 15647.32 | 144.16 | 13326.43 | 141.82 | 0 | 0 |
| C206 | 20991.54 | 0.334 | 16591.19 | 0.265 | 20991.54 | 146.94 | 16591.19 | 144.14 | 0 | 0 |
| C207 | 17919.79 | 0.290 | 13835.00 | 0.25 | 17919.79 | 156.94 | 13835.00 | 156.62 | 0 | 0 |
| C208 | 13648.89 | 0.268 | 13326.53 | 0.25 | 13648.89 | 143 | 13326.53 | 142.03 | 0 | 0 |
| **Average** | **18398.97** | **0.317** | **15757.98** | **0.263** | **18398.97** | **152.32** | **15757.98** | **150.24** | **0** | **0** |
| R201 | 36043.16 | 0.426 | 36043.1 | 0.234 | 36043.1 | 178.57 | 36043.1 | 176.51 | 0 | 0 |
| R202 | 36043.16 | 0.273 | 36043.16 | 0.187 | 36043.16 | 180.42 | 36043.16 | 178.42 | 0 | 0 |
| R203 | 31809.51 | 0.578 | 31809.51 | 0.5 | 31809.51 | 182.57 | 31809.51 | 180.68 | 0 | 0 |
| R204 | 28968.51 | 0.729 | 28968.51 | 0.671 | 28968.51 | 181.73 | 28968.51 | 179.92 | 0 | 0 |
| R205 | 22839.41 | 0.376 | 22839.41 | 0.312 | 22839.41 | 178.57 | 22839.41 | 176.26 | 0 | 0 |
| R206 | 27890.90 | 0.3 | 27890.90 | 0.265 | 27890.90 | 142.94 | 27890.90 | 139.31 | 0 | 0 |
| R207 | 26693.56 | 0.376 | 26693.56 | 0.343 | 26693.56 | 146.7 | 26693.56 | 145.43 | 0 | 0 |
| R208 | 15476.88 | 0.689 | 13198.77 | 0.484 | 15476.88 | 129.81 | 13198.77 | 128.87 | 0 | 0 |
| R209 | 29602.88 | 0.251 | 29602.88 | 0.218 | 29602.88 | 130.02 | 29602.88 | 129.60 | 0 | 0 |
| R210 | 26091.07 | 0.381 | 26091.07 | 0.343 | 26091.07 | 130.19 | 26091.07 | 128.59 | 0 | 0 |
| R211 | 29954.17 | 0.328 | 29954.17 | 0.296 | 29954.17 | 126.85 | 29954.17 | 125.71 | 0 | 0 |
| **Average** | **28310.29** | **0.427** | **28103.18** | **0.350** | **28310.28** | **155.30** | **28103.18** | **153.57** | **0** | **0** |
| RC201 | 39036.03 | 0.203 | 39036.03 | 0.187 | 39036.03 | 133.40 | 39036.03 | 132.20 | 0 | 0 |
| RC202 | 27261.97 | 0.221 | 27261.97 | 0.203 | 27261.97 | 135.65 | 27261.97 | 134.81 | 0 | 0 |
| RC203 | 27317.06 | 0.484 | 27317.06 | 0.453 | 27317.06 | 127.17 | 27317.06 | 126.59 | 0 | 0 |
| RC204 | 25032.85 | 0.371 | 25032.85 | 0.343 | 25032.85 | 135.14 | 25032.85 | 134.10 | 0 | 0 |
| RC205 | 36512.4 | 0.215 | 36512.4 | 0.171 | 36512.4 | 132.46 | 36512.4 | 131.28 | 0 | 0 |
| RC206 | 27428.40 | 0.353 | 27428.40 | 0.328 | 27428.40 | 135.14 | 27428.40 | 133.96 | 0 | 0 |
| RC207 | 32745.40 | 0.343 | 32745.40 | 0.312 | 32745.40 | 135.01 | 32745.40 | 134.35 | 0 | 0 |
| RC208 | 33066.70 | 0.431 | 33066.70 | 0.390 | 33066.70 | 132.79 | 33066.70 | 132.15 | 0 | 0 |
| **Average** | **31050.10** | **0.327** | **31050.10** | **0.298** | **31050.10** | **133.34** | **31050.10** | **132.43** | **0** | **0** |
| **Avg All** | **25919.79** | **0.357** | **24970.42** | **0.304** | **25919.78** | **146.99** | **24970.42** | **145.41** | **0** | **0** |

Table 5 and Table 6 report the computational results for the modified type 1 (R1, C1, RC1) and type 2 (R2, C2, RC2) instances with $\Delta = 30\%, 50\%, 70\%, 90\%$, respectively. We have chosen the best execution for each instance and grouped them to calculate the average values of each group with different degrees of dynamism. C2 is the dynamic degree of the system. From C3 to C10 are the computational result and they are the average objective value of the solutions found in 10 runs (Obj_A), the average computational time in seconds required to find the solutions in 10 runs (CPU_A), the best objective value of the solutions found in 10 runs (Obj_B), and the best computational time in seconds required to find the solutions in 10 runs (CPU_B), respectively. C11 is the relative error between C5 and C9 (i.e. C11=(C5-C9)/C9*100%), and C12 is the one between C3 and C7 (i.e. C12=(C3-C7)/C7*100%).

Table 5. Comparison results of type 1 instances for other dynamic degrees.

| Type (1) | $\Delta$ (2) | EIADR | | | | EIADR +TS | | | | BRE (11) | ARE (12) |
|---|---|---|---|---|---|---|---|---|---|---|---|
| | | Obj_A (3) | CPU_A (4) | Obj_B (5) | CPU_B (6) | Obj_A (7) | CPU_A (8) | Obj_B (9) | CPU_B (10) | | |
| R1 | 30% | 96472.49 | 0.532 | 93890.0 | 0.442 | 96472.49 | 152.78 | 93890.0 | 151.14 | 0 | 0 |
| | 50% | 82957.83 | 0.827 | 82196.23 | 0.75 | 82957.83 | 151.45 | 82196.23 | 149.73 | 0 | 0 |
| | 70% | 77731.85 | 0.938 | 75202.94 | 0.858 | 77731.85 | 150.03 | 75202.94 | 147.79 | 0 | 0 |
| | 90% | 159630.5 | 1.075 | 159228.7 | 1.005 | 150190.2 | 152.56 | 146584.7 | 150.80 | 8.625 | 6.285 |
| C1 | 30% | 68085.37 | 0.665 | 63837.09 | 0.513 | 67615.02 | 152.00 | 63520.24 | 149.15 | 0.498 | 0.695 |
| | 50% | 57473.44 | 1.034 | 54711.88 | 0.829 | 56979.92 | 150.24 | 54220.93 | 144.89 | 0.905 | 0.866 |
| | 70% | 52106.59 | 1.102 | 49553.77 | 0.907 | 51297.78 | 147.83 | 48810.52 | 141.57 | 1.522 | 1.576 |
| | 90% | 46828.82 | 1.268 | 46440.27 | 1.006 | 46039.07 | 149.67 | 45540.78 | 145.02 | 1.975 | 1.715 |
| RC1 | 30% | 74523.11 | 0.499 | 71717.89 | 0.464 | 74523.11 | 143.38 | 71717.89 | 140.92 | 0 | 0 |
| | 50% | 72523.85 | 0.838 | 71317.23 | 0.744 | 72523.85 | 151.98 | 71317.23 | 150.05 | 0 | 0 |
| | 70% | 75217.93 | 0.937 | 73261.84 | 0.884 | 75217.93 | 151.25 | 73261.84 | 150.10 | 0 | 0 |
| | 90% | 75225.42 | 1.057 | 75187.65 | 1.009 | 75225.42 | 152.39 | 75187.65 | 150.67 | 0 | 0 |
| **Ave** | | **78231.43** | **0.897** | **76378.79** | **0.784** | **77231.20** | **150.46** | **75120.91** | **147.65** | **1.674** | **1.295** |

Table 6. Comparison results of type 2 instances for other dynamic degrees.

| Type (1) | Δ (2) | EIADR | | | | EIADR +TS | | | | BRE (11) | ARE (12) |
|---|---|---|---|---|---|---|---|---|---|---|---|
| | | Obj_A (3) | CPU_A (4) | Obj_B (5) | CPU_B (6) | Obj_A (7) | CPU_A (8) | Obj_B (9) | CPU_B (10) | | |
| R2 | 30% | 21091.56 | 1.049 | 19190.12 | 0.830 | 21091.56 | 154.47 | 19190.12 | 152.31 | 0 | 0 |
| | 50% | 18542.97 | 1.345 | 18542.97 | 1.235 | 18542.97 | 155.07 | 18542.97 | 153.26 | 0 | 0 |
| | 70% | 16824.98 | 2.085 | 15782.98 | 1.909 | 16824.98 | 153.86 | 15782.98 | 151.47 | 0 | 0 |
| | 90% | 17592.60 | 2.178 | 17305.36 | 1.985 | 17592.60 | 152.91 | 17305.36 | 150.81 | 0 | 0 |
| C2 | 30% | 16639.19 | 0.812 | 14996.68 | 0.675 | 16639.19 | 150.25 | 14996.68 | 148.89 | 0 | 0 |
| | 50% | 15664.51 | 1.468 | 14347.97 | 1.226 | 15664.51 | 149.89 | 14347.97 | 147.65 | 0 | 0 |
| | 70% | 15070.10 | 1.072 | 14882.28 | 0.998 | 15070.10 | 145.33 | 14882.28 | 143.03 | 0 | 0 |
| | 90% | 15155.28 | 1.157 | 14529.24 | 1.050 | 15155.28 | 141.14 | 14529.24 | 139.34 | 0 | 0 |
| RC2 | 30% | 25597.84 | 0.787 | 25469.95 | 0.703 | 25597.84 | 136.77 | 25469.95 | 135.97 | 0 | 0 |
| | 50% | 22215.23 | 1.311 | 22066.07 | 1.216 | 22215.23 | 143.74 | 22066.07 | 141.44 | 0 | 0 |
| | 70% | 21216.35 | 1.501 | 20822.01 | 1.416 | 21216.35 | 148.14 | 20822.01 | 146.80 | 0 | 0 |
| | 90% | 20925.52 | 1.689 | 20385.09 | 1.609 | 20925.52 | 150.74 | 20385.09 | 149.50 | 0 | 0 |
| **Ave** | | **18878.01** | **1.371** | **18193.39** | **1.237** | **18878.01** | **148.52** | **18193.39** | **146.70** | **0** | **0** |

Based on Tables 3-6, we can draw the following conclusions.

(i) Compared with the EIADR+TS method, the proposed EIADR method in this study exhibits excellent properties for type 2 instances, while the solution performance is slightly worse for type 1 with a smaller vehicle capacity. As shown in Table 3, the relative errors (BRE and ARE) at 10% dynamic degree are 0.256 and 0.189, respectively. For the remaining dynamic rates, the relative errors rise to BRE=1.674 and ARE=1.295 (see Table 5), while the relative error for type 2 instances is 0%. This may be due to the narrow time window of type 1 instances, which makes it easier to miss the optimal solution when planning dynamic requests using the EIADR method. Therefore, this method is more effective in solving dynamic problems with wider time windows.

(ii) The effectiveness of using the EIADR method for planning type C instances with node cluster distribution is not effective compared to R instances and RC instances. Moreover, compared to the results of using the TS algorithm for further optimization, the relative error can reach BRE=16.74 and ARE=9.042, as shown in Table 3.

(iii) One can find that as the dynamic degree of the system increases, the optimization effect of the instance also decreases, as is notably observed in the results of C1 instances in Table 5. When the R1 type data reaches a dynamic degree of 90%, the relative error BRE reaches 8.625%, as shown in Table 5. The efficiency of the EIADR method proposed in this paper in solving DAMRRP-TWRP largely depends on dynamic degree of problems.

(iv) The results in Table 2-5 show that our EIADR method is very fast and robust in solving the DAMRRP-TWRP, as the total average computational time per instance is about 1 second which is much less than that of EIADR+TS method and the best objective values is very close to the average objective values. In order words, when a patient makes a request, he/she can get his/her service order after 1 second. Therefore, for situations with high timeliness requirements, we can choose the EIADR method to directly handle dynamically arriving requests.

*5.3.2. Computational results for instances*

In this subsection, some analysis is conducted using RC202_10% to validate the effectiveness of the proposed algorithm. The results are presented in Table 7. Firstly, the TS algorithm is applied to plan for all known static requests, generating a set of prior routes, such as Route 1 in Table 7, where 0 represents depot and 1-90 denotes all static requests. Subsequently, the ten dynamic requests, 92-101, appear sequentially during the operation of AMRs and are marked with an asterisk (*). Based on the prior routes, the EIADR method is employed to optimize the received dynamic requests in real-time, resulting in the final scheduling scheme, i.e., Route 2 in Table 7, where dynamic request 99 is rejected due to excessively high delay penalty cost. These results show that the EIADR method can solve the real-time requirements of the DAMRRP-TWRP.

Table 7. Case analysis for RC202_10%.

| No. | Route 1 (static request) | Route 2 (including dynamic request) |
| --- | --- | --- |
| AMR1 | 0-2-1-3-5-7-12-21-25-30-33-43-45-47-54-65-68-74-75-77-83-85-86-87-88-0 | 0-2-1-3-5-7-12-21-25-30-33-43-45-47-54-65-68-74-75-77-83-85-86-87-88-0 |
| AMR2 | 0-59-14-57-32-26-27-56-58-64-55-36-71-6-41-10-20-51-60-50-4-17-70-53-24-0 | 0-59-14-57-32-26-27-56-58-64-55-36-**97**$^*$-71-6-41-10-20-51-60-50-4-17-70-53-24-0 |
| AMR3 | 0-35-37-42-62-15-11-19-61-46-18-22-76-44-39-49-67-23-81-84-63-0 | 0-35-37-42-62-15-11-19-61-46-18-22-76-44-39-49-**100**$^*$-67-23-81-84-63-0 |
| AMR4 | 0-28-40-29-78-79-9-52-31-34-90-72-0 | 0-**92**$^*$-**93**$^*$-28-40-29-78-79-9-52-31-34-**96**$^*$-90-72-0 |
| AMR5 | 0-16-80-69-89-73-38-13-0 | 0-**95**$^*$-16-80-69-89-73-38-**98**$^*$-13-**101**$^*$-0 |
| AMR6 | 0-66-82-8-0-48-0 | 0-**94**$^*$-66-82-8-0-48-0 |

### 5.3.3. Management insights

In this subsection, we discuss management insights with respect to the AMR's safety stock settings. To better describe the impact of different dynamic degrees on the quality of dynamic requests of AMRs services, we introduce service rate $\tau = \frac{\sum_{i \in R_3} y_i}{n_d} *100\%$, where $\sum_{i \in R_3} y_i$ is the total number of dynamic requests served by AMRs.

For the impact of different safety stocks on service rates, we analyzed instances of type 1 (R1, C1, RC1) and type 2 (R2, C2, RC2), respectively. The results are shown in Figure 2. Figure 2 (a) and (b) show the impact of safety stock on service rates for type 1 and type 2 instances under five dynamic degrees, respectively.

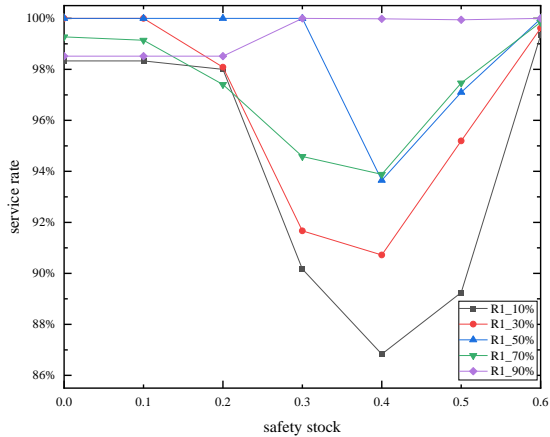

(1) instance R1

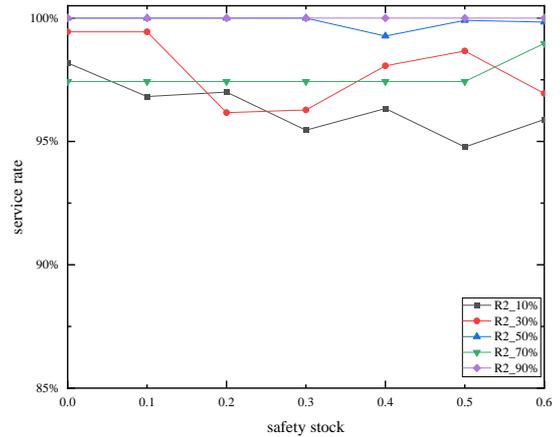

(4) instance R2

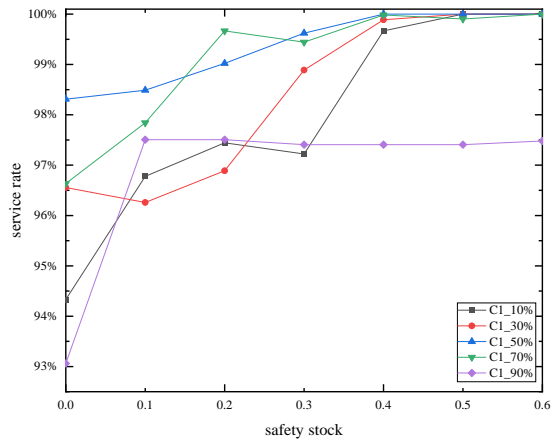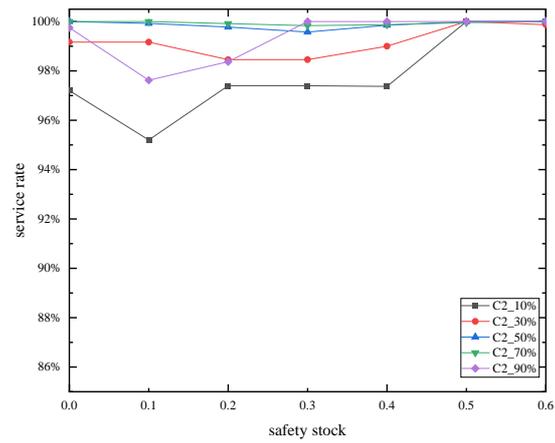

(2) instance C1          (5) instance C2

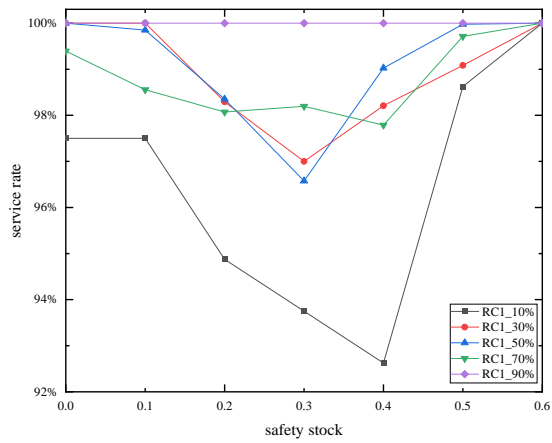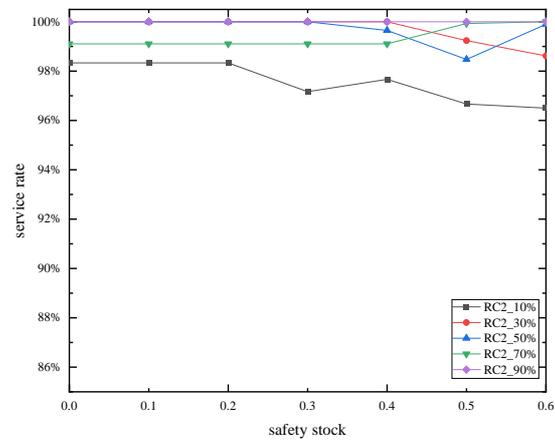

(3) instance RC1          (6) instance RC2

(a) Type 1 instances          (b) Type 2 instances

Figure 2. The impact of different replenishment strategies on service rates.

Figure 2 (a) (1)-(3) shows the service rate of the type C1 instance significantly improves with an increase in safety stock, regardless of whether the dynamic rate is 10% or 90%. This may be because the locations of the nodes in type C1 problems are clustered distribution and the capacity of the vehicle is relatively small, which necessitates an increase in the safety stock of AMR to ensure efficient processing of dynamic requests and improve the service rate of requests. However, for type R1 and RC1 instances with similarly small vehicle capacity, there exists a threshold when the safety stock is approximately 0.3 to 0.4, which leads to the minimum service rate. Therefore, for type

R1 and RC1 with non-clustered node distribution, when planning the service route for static requests, the safety stock of AMRs should be set to the minimum, that is, to make each AMR fully serve static requests as much as possible, so that more dynamic requests are scheduled to the new AMR to prevent excessive time-window violations. Alternatively, the safety stock of AMRs can be set higher, which means more AMRs are deployed to handle arriving dynamic requests, thereby improving the service rate. Therefore, for the case of small capacity and relatively clustered node distribution, increasing the safety stock of AMR is an effective way to improve the service rate of medical requests.

As can be seen from Figure 2 (b) (4)-(6), the service rate of type 2 instances is not significantly affected by safety stock. This may be because type 2 instances have a relatively larger vehicle capacity and a wider time window compared to type 1 instances, resulting in a smaller impact of the safety stock setting on the arrival of dynamic requests, making the service rate relatively unaffected by the safety stock. Therefore, for dynamic problems with relatively large vehicle capacity and a wider time window, the safety stock of AMR can be set based on the minimum cost of introducing AMR by hospitals.

## 6. Conclusion

This paper studied the DAMRRP-TWRR based on a real-life application, which takes into account three dynamic factors: stochastic service time, stochastic travel time, and dynamic requests. The DAMRRP-TWRR opens the door to the scheduling problem of AMRs participating in medical services. In the paper, we have established a two-stage mathematical programming model for DARRP-TWRR and propose some service strategies for scheduling dynamic request services for DAMRRP-TWRR. A fast and efficient EIADR method is proposed based on these service strategies. Compared with the EIADR+TS algorithm, which optimizes the route using the tabu search algorithm after using the EIADR method, the computational results show that the EIADR method is very fast and robust in solving the DAMRRP-TWRP. For situations with high timeliness requirements, we can choose the EIADR method to directly handle dynamically arriving requests. Finally, we provide managerial insights concerning the AMR's safety stock settings in situations where the case of small vehicle capacity and relatively clustered node distribution. We highlighted the fact that to improve the service rate of dynamic requests, AMR capacity space should be appropriately reserved to better handle randomly arrived medical requests.

This work only considers the distribution of required drugs by AMRs in hospitals from one delivery point (such as a pharmacy) to various request points (such as wards). It could happen that AMRs need to depart from multiple depots and perform simultaneous pickup and delivery transportation. Therefore, future research can explore the dynamic pickup and delivery AMR scheduling problem with multi-depot.

**Disclosure statement**

# Appendix A.

Table A1. The pseudocode of Greedy Insertion Algorithm for a Feasible Solution.

**Algorithm A.1. Greedy Insertion Algorithm for a Feasible Solution**

**Step 1** Arrange the lower bound $e_i$ of the time window of all requests in ascending order to generate a sequence of service priorities for requests $P = \{p_1, p_2, \ldots, p_i, \ldots, p_n\}$. Let $i = 1$ and the route set is empty ($R_{current} = \Phi$).

**Step 2** Assign $p_1$ to a new AMR, and generate sub-route $r_1 = \{0 \to p_1\}$.
$R_{current} \leftarrow R_{current} \cup r_1$.

**Step 3** If $i \leq n$
    then repeat **Step 4**
  Else
    then go to **Step 5**.

**Step 4** If there exists a sub-route $r_j \in R_{current}, j = 1, 2, \ldots$ such that the request $p_i$ can be served at the end of $r_j$, and the constraints of the time window and capacity are satisfied.
    then update sub-route $r_j = \{r_j \to p_i\}$, set $i \leftarrow i + 1$.
  Else
    then assign $p_i$ to a new AMR and generate a new sub-route $r_{new} = \{0 \to p_i\}$.
    Set $i \leftarrow i + 1$, $R_{current} \leftarrow R_{current} \cup r_{new}$.
  End

**Step 5** Let $R_{current}$ be the initial feasible solution $x_0$, and output the corresponding objective function $F(x_0)$.

Table A2. The pseudocode of TS for the AMR Scheduling Problem.

**Algorithm A.2. Tabu Search (TS) for the AMR Scheduling Problem**

**Step 1** Initialize the weights of the three operators $\rho_i = 1, i = 1, 2, 3$ and $\delta_1, \delta_2$, three tabu lists $B_i, i = 1, 2, 3$, tabu tenure $t$ and the algorithm iterations $ite = 1$.

**Step 2** Let $x_0$ (obtained from **Algorithm A.1**) be the initial solution of the algorithm, i.e., $x_{current} = x_0$, $F(x_{current}) = F(x_0)$. And the present best solution is $x_{best} = x_0$, $F(x_{best}) = F(x_0)$.

**Step 3** Based on the present solution $x_{current}$, use the roulette wheel procedure for the three operators (swap*, 2-opt*, and relocation*) to select the operator $i$ and generate the corresponding neighborhood $N(x_{current})$. The depot insertion operator is used to repair neighborhood solutions. Find the optimal objective value $F(x_{c\_best})$ and the corresponding neighborhood action $b_i(j_1^*, j_2^*)$, where $x_{c\_best} \in N(x_{current})$.

**Step 4** if $F(x_{c\_best}) < F(x_{best})$

**then** $F(x_{best}) = F(x_{c\_best})$, $x_{best} = x_{c\_best}$, $x_{current} = x_{c\_best}$, $F(x_{current}) = F(x_{c\_best})$, the weight $\rho_i$ of the selected operator $i$ is updated to $\rho_i = \rho_i + \delta_1$ through the reward scheme, and update the tabu list. Let

$$B_i(j_1^*, j_2^*) = \begin{cases} 0, \text{ if } B_i(j_1^*, j_2^*) \neq 0 \\ t, \text{ if } B_i(j_1^*, j_2^*) = 0 \end{cases}, B_i(j_1, j_2) = B_i(j_1, j_2) - 1 \text{ if}$$

$B_i(j_1, j_2) \neq 0$ and $j_1 = j_1^*, j_2 \neq j_2^*$.

**Elseif** $F(x_{c\_best}) \geq F(x_{best})$

**then** find the optimal solution $x'_{c\_best}$ generated by the neighborhood action $b_i(j_1', j_2')$ for $\forall B_i(j_1', j_2') = 0$. Let $x_{current} = x'_{c\_best}$, $F(x_{current}) = F(x'_{c\_best})$, the weight $\rho_i$ of the selected operator $i$ is updated to $\rho_i = \rho_i + \delta_2$ through the reward scheme, and update the tabu list. Let $B_i(j_1', j_2') = t$, $B_i(j_1, j_2) = B_i(j_1, j_2) - 1$ if $B_i(j_1, j_2) \neq 0$.

**End**

**Step 5** Let $ite = ite + 1$. The probability of each operator $P_i = \dfrac{\rho_i}{\sum\limits_{j=1}^{3} \rho_j}$ is updated after every 10 iterations. Repeat **Steps 3-4**. When $ite = N$, stop the algorithm and output the present optimal solution $F(x_{best})$ and the corresponding set of routes $x_{best}$.